\documentclass[graybox]{svmult}

\usepackage{mathptmx}       
\usepackage{helvet}         
\usepackage{courier}        
\usepackage{makeidx}         
\usepackage{graphicx}        
\usepackage{multicol}        
\usepackage[bottom]{footmisc}

\usepackage{cite}
\usepackage{fancyvrb}        
\usepackage{tabularx}
\usepackage{siunitx}
\usepackage{multirow}
\usepackage{amsmath,amssymb} 
\newtheorem{rem}{Remark}

\begin{document}

\title*{Towards a Parallel-in-Time Calculation of Time-Periodic Solutions with Unknown Period}
\titlerunning{Parallel-in-Time Calculation of Time-Periodic Solutions with Unknown Period}
\author{Iryna Kulchytska-Ruchka 
\and
Sebastian Sch{\"o}ps 
}
\institute{Iryna Kulchytska-Ruchka, Sebastian Sch{\"o}ps \at Computational Electromagnetics Group, Technical University of Darmstadt, Schlossgartenstr. 8, 64289 Darmstadt, Germany,  
\email{\{kulchytska, schoeps\}@temf.tu-darmstadt.de}
 }

%
%
\maketitle

\abstract*{}

\abstract{This paper presents a novel parallel-in-time algorithm able to compute time-periodic solutions 
of problems where the period is not given. Exploiting the idea of the multiple shooting 
method, the proposed approach calculates the initial values at each subinterval as well as the corresponding 
period iteratively. As in the Parareal method, parallelization in the time domain is performed using discretization on a two-level grid. A special linearization of the time-periodic system on the coarse grid is introduced {to speed up the computations}. The iterative algorithm is verified via its application to the Colpitts oscillator model. 
}

\section{Introduction}
\label{kulchytska-ruchka:sec_general}
Steady-state analysis is {a common task} in electrical engineering, 
{for example,} during the initial design stages of, 
e.g., electric {circuits or} motors. Classical sequential time 
stepping {may} lead to lengthy transient computation particularly 
when the underlying dynamical system possesses a large time constant. 
Various approaches for efficient steady-state calculation are known from the literature. For instance, {clever methods to choose} the starting value \cite{Bermudez_2019aa} or an explicit error correction \cite{Takahashi_2010aa} could accelerate the time-domain calculation considerably.

A powerful tool for speeding up the classical time stepping is the class of parallel-in-time methods, such as the multigrid reduction in time \cite{Falgout_2014aa} or Parareal \cite{Lions_2001aa}. Originating from the multiple shooting method \cite{Morrison_1962aa}, they are based on the splitting of the considered time interval into several windows and updating the solution at synchronization points iteratively. 
The use of coarse and fine discretizations {propagates quickly low-frequency } information {of} the solution {using} a cheap sequential solver followed by a very accurate result with a  precise fine solver applied in parallel.

Another direction of obtaining the steady state is based on the
solution of the joint space- and time-discrete time-periodic system 
formulated on the whole period \cite{Hara_1985aa}. There the initial and final values are coupled through the prescribed periodicity condition. An obstacle within the solution of the periodic problem in the time domain 
becomes the large size of the system matrix as well as its special block structure due to the interdependence of the solution vectors over the period. 
To deal with this difficulty a frequency domain approach was proposed in 
\cite{Biro_2006aa}. In case of linear problems, the method takes advantage of the block-cyclic matrix
structure by applying the discrete Fourier transform. 
It fully decouples the variables, thereby allowing 
for the separate solution {of} each harmonic coefficient. 
This approach was further extended and incorporated into the 
Parareal framework by the authors in \cite{Kulchytska-Ruchka_2019ag}. There, a simplified Newton-based iterative algorithm was presented together with its convergence analysis for the efficient treatment of nonlinear problems.

Solutions of time-periodic problems become much more challenging when the 
period is not given. Such situation occurs, e.g., when dealing with an autonomous system \cite{Deuflhard_1984aa}. In contrast to a non-autonomous problem, the periodicity cannot be determined from the applied excitation. 
This paper proposes a numerical algorithm capable of determining an appropriate period automatically using parallelization in the time domain. Extending the idea of the multiple shooting method we include the unknown period together with multiple initial 
values as the sought parameters into the iterative procedure. Verification of the presented approach is illustrated through its application to the Colpitts oscillator model {\cite{Kampowsky_1992aa}}. 

The paper is organized as follows. Section~\ref{kulchytska-ruchka:multishoot} describes the basis of the multiple shooting approach including the unknown period as an additional variable. This is further expanded to the family of the Parareal-based methods in Section~\ref{kulchytska-ruchka:parallel-in-time}. Section~\ref{kulchytska-ruchka:numerics} applies the proposed parallel-in-time approach to the Colpitts oscillator model using a particular linerization on the coarse level. The paper is finally summarized in Section~\ref{kulchytska-ruchka:conclusions}.

\section{Multiple shooting with unknown period}
\label{kulchytska-ruchka:multishoot}
We consider the following time-periodic problem for a system of ordinary differential equations (ODEs)
\begin{equation}
\label{kulchytska-ruchka:TP_pbm}
\begin{aligned}
\mathbf{M}\tilde{\mathbf{u}}^{\prime}(t)&=\mathbf{f}(\tilde{\mathbf{u}}(t)),\ \ t\in(0,T)\\
\tilde{\mathbf{u}}(0)&=\tilde{\mathbf{u}}(T),
\end{aligned}
\end{equation}
where the period $T>0$ and the vector $\tilde{\mathbf{u}}:[0,T]\to\mathbb{R}^d,$ $d\geq1$ are sought. {$\mathbf{M}$ is a given {non-singular} mass matrix, $\mathbf{f}$ is a bounded and Lipschitz continous right-hand side (RHS) function.}  
Following \cite{Deuflhard_1984aa} we incorporate the period $T$ as an unknown parameter by performing the change of variables
\begin{equation}
\label{kulchytska-ruchka:scaling}
[0,T]\ni t\mapsto\tau:=t/T\in[0,1].
\end{equation}
The problem \eqref{kulchytska-ruchka:TP_pbm} is thereby transformed into the equivalent one: find $T>0$ and ${\mathbf{u}}:[0,1]\to\mathbb{R}^d$ such that 
\begin{equation}
\label{kulchytska-ruchka:TP_pbm_unitInt}
\begin{aligned}
\mathbf{M}{\mathbf{u}}^{\prime}(\tau)&=T\mathbf{f}({\mathbf{u}}(\tau)),\ \ \tau\in(0,1)\\
{\mathbf{u}}(0)&={\mathbf{u}}(1).
\end{aligned}
\end{equation}
The unit interval $[0,1]$ is then partitioned into $N$ windows by the nodes $0=\tau_0<\tau_1<\dots<\tau_N=1.$ The $n$-th subinterval has length $\Delta\tau_n=\tau_n-\tau_{n-1},$ for $n=1,\dots,N.$

\begin{figure}[t]
\centering
\includegraphics{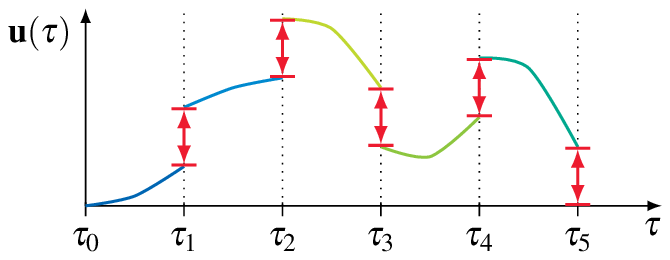}
\caption{Example of the interval splitting within the multiple shooting for $N=5$. The mismatches at the {synchronization points} $\tau_n,$ $n=1,\dots,N-1$ together with the periodicity jump between the solution at $\tau_0$ and $\tau_N$ are eliminated (up to a prescribed tolerance) by solving the root-finding problem.}
\label{kulchytska-ruchka:fig_1}   
\end{figure}

For a given discrete variable $\mathbf{U}_{n-1}$, we consider an initial value problem (IVP) on the window $(\tau_{n-1},\tau_n]$
\begin{equation}
\label{kulchytska-ruchka:IVP_window}
\begin{aligned}
\mathbf{M}{\mathbf{u}}^{\prime}_n(\tau)&=T\mathbf{f}({\mathbf{u}_n}(\tau)),\ \ \tau\in(\tau_{n-1},\tau_n]\\
{\mathbf{u}_n}(\tau_{n-1})&={\mathbf{U}_{n-1}}
\end{aligned}
\end{equation}
and let $\mathcal{F}(\tau_{n},\tau_{n-1},\mathbf{U}_{n-1},T)$ denote the solution operator of \eqref{kulchytska-ruchka:IVP_window} for $n=1,\dots,N.$ A sketch of the piecewise-defined solution due to the interval splitting is shown in Fig.~\ref{kulchytska-ruchka:fig_1}. In order to eliminate the jumps at the synchronization points $\tau_n,$ $n=1,\dots,N-1$ as well as the difference between the initial value at $\tau_0$ and the final one at $\tau_N$ the matching conditions:
\begin{equation}\label{kulchytska-ruchka:jumps}
\varPhi(\mathbf{z}):=
 \left\{
\begin{aligned}
\mathcal{F}(\tau_{N},\tau_{N-1},\mathbf{U}_{N-1},T)-\mathbf{U}_{0}&=0,\\
\mathcal{F}(\tau_{n},\tau_{n-1},\mathbf{U}_{n-1},T)-\mathbf{U}_{n}&=0,\ \ n=1,\dots,N-1
\end{aligned} \right.
\end{equation}
{have to be satisfied,} where $\mathbf{z}=\left[\mathbf{U}_0^{\top},\dots,\mathbf{U}_{N-1}^{\top},T\right]^{\!\top}\!.$ System \eqref{kulchytska-ruchka:jumps} represents the root-finding problem for the mapping $\varPhi:\mathbb{R}^{Nd+1}\rightarrow\mathbb{R}^{Nd}.$ The Jacobian 
of $\varPhi$ is given by
\begin{equation}
\mathbf{J}_{\varPhi}(\mathbf{z})=
\begin{bmatrix}
-\mathbf{I} & & & \mathbf{G}_{N} & \mathbf{g}_{N}\\
\mathbf{G}_{1} & -\mathbf{I} & & & \mathbf{g}_{1}\\
& \ddots & \ddots &  & \vdots\\
& & \mathbf{G}_{N-1} & -\mathbf{I} & \mathbf{g}_{N-1}
\end{bmatrix},
\end{equation}
where
\begin{equation}
\label{kulchytska-ruchka:Jacobian_Phi}
\mathbf{G}_n=\frac{\partial\mathcal{F}}{\partial\mathbf{U}_{n-1}}(\tau_{n},\tau_{n-1},\mathbf{U}_{n-1},T),\qquad 
\mathbf{g}_n=\frac{\partial\mathcal{F}}{\partial T}(\tau_{n},\tau_{n-1},\mathbf{U}_{n-1},T),
\end{equation}
$n=1,\dots,N$ and $\mathbf{I}$ is the identity matrix. The root $\mathbf{z}$ of \eqref{kulchytska-ruchka:jumps} can then be calculated using the Newton method, i.e., for a given $\mathbf{z}^{(0)}$ and $k=0,1,\dots$ solution $\mathbf{z}^{(k+1)}$ at the iteration $k+1$ is updated through
\begin{align}
\mathbf{J}_{\varPhi}\bigl(\mathbf{z}^{(k)}\bigr)\Delta\mathbf{z}^{(k)}&=-{\varPhi}\bigl(\mathbf{z}^{(k)}\bigr),
\label{kulchytska-ruchka:Newton_multiShoot}\\
\mathbf{z}^{(k+1)}&=\Delta\mathbf{z}^{(k)}+\mathbf{z}^{(k)}.
\label{kulchytska-ruchka:Newton_multiShoot_update}
\end{align}
Note that due to the introduction of the additional variable $T$ the system of equations \eqref{kulchytska-ruchka:Newton_multiShoot} is underdetermined and can be solved, e.g., by calculating the Moore-Penrose pseudoinverse. {A generalized eigenvalue-based gauging as well as the corresponding theory for the Moore-Penrose pseudoinversion was presented in \cite{Boroujeni_2010aa}.} We note that in case when the size of \eqref{kulchytska-ruchka:Newton_multiShoot} is large it can be condensed to a $d$-dimensional system with $d+1$ unknowns by block Gaussian elimination \cite{Deuflhard_1984aa}.

\section{Periodic time-parallelization with coarse grid correction}
\label{kulchytska-ruchka:parallel-in-time}
Inheriting the idea of the Parareal algorithm \cite{Lions_2001aa,Gander_2007ab} we approximate the derivative in $\mathbf{G}_n$ \eqref{kulchytska-ruchka:Jacobian_Phi} in a finite difference way using a coarse propagator $\mathcal{G},$ i.e., for the iteration $k$ and $n=1,\dots,N$
\begin{equation}
\label{kulchytska-ruchka:parareal_approxG}
\begin{aligned}
\mathbf{G}_n^{(k)}\Delta\mathbf{U}_{n-1}^{(k)}&=\frac{\partial\mathcal{F}}{\partial\mathbf{U}_{n-1}}\bigl(\tau_{n},\tau_{n-1},\mathbf{U}_{n-1}^{(k)},T\bigr)\bigl[\mathbf{U}_{n-1}^{(k+1)}-\mathbf{U}_{n-1}^{(k)}\bigr]\\
&\approx\mathcal{G}\bigl(\tau_{n},\tau_{n-1},\mathbf{U}_{n-1}^{(k+1)},T\bigr)-\mathcal{G}\bigl(\tau_{n},\tau_{n-1},\mathbf{U}_{n-1}^{(k)},T\bigr).
\end{aligned}
\end{equation}
{Similar to} the fine propagator $\mathcal{F},$ the operator $\mathcal{G}$ solves the IVP \eqref{kulchytska-ruchka:IVP_window} on each time window. However, in contrast to the fine solver the coarse propagator has a considerably lower precision, e.g., it uses a lower-order time integrator or bigger time steps. Substituting \eqref{kulchytska-ruchka:parareal_approxG} into \eqref{kulchytska-ruchka:Newton_multiShoot} we obtain the periodic Parareal with periodic coarse problem \cite{Gander_2013ab} with unknown period {(PP-PC-UP)}:
\begin{equation}
\label{kulchytska-ruchka:PPPC_unknownT}
\setlength{\arraycolsep}{2.5pt}
\begin{bmatrix}
-\mathbf{I} &    & & {\mathcal{G}_N\left(\cdot,T\right)} & \mathbf{g}^{(k)}_{N} \\
{\mathcal{G}_1\left(\cdot,T\right)}  & -\mathbf{I} &  &  & \mathbf{g}^{(k)}_{1}\\
 &  \ddots   &  \ddots & &\vdots \\
 &      &  {\mathcal{G}_{N-1}\left(\cdot,T\right)} & -\mathbf{I} & \mathbf{g}^{(k)}_{N-1}
\end{bmatrix}
	\begin{bmatrix}
	\mathbf{U}_0^{(k+1)}\\
	\mathbf{U}_1^{(k+1)}\\
	\vdots\\
	\mathbf{U}_{N-1}^{(k+1)}\\
	T^{(k+1)}
	\end{bmatrix}=\begin{bmatrix}
		\mathbf{b}_N^{(k)}\\
		\mathbf{b}_1^{(k)}\\
		\vdots\\
		\mathbf{b}_{N-1}^{(k)}
\end{bmatrix},
\end{equation}
where we denote $\mathbf{b}_n^{(k)}:=\mathbf{g}_n^{(k)}T^{(k)}+\mathcal{G}\bigl(\tau_n, \tau_{n-1},\mathbf{U}^{(k)}_{n-1},T\bigr)-\mathcal{F}\bigl(\tau_n, \tau_{n-1},\mathbf{U}^{(k)}_{n-1},T\bigr)$ and {$\mathcal{G}_n\bigl(\cdot,T\bigr):=\mathcal{G}\bigl(\tau_n,\tau_{n-1},\cdot,T\bigr)$} for $n=1,\dots,N.$ Following \cite{Deuflhard_1984aa} we have
\begin{equation}
\label{kulchytska-ruchka:approx_g}
\begin{aligned}
\mathbf{g}_n^{(k)}&=\frac{\partial\mathcal{F}}{\partial T}\bigl(\tau_{n},\tau_{n-1},\mathbf{U}_{n-1}^{(k)},T\bigr)=\frac{\partial}{\partial T}\Bigl[\mathbf{U}_{n-1}^{(k)}+\mathbf{M}^{-1}\int\limits_{\tau_{n-1}}^{\tau_{n}}T\mathbf{f}({\mathbf{u}}(\tau))d\tau\Bigr]\\
&=\mathbf{M}^{-1}\int\limits_{\tau_{n-1}}^{\tau_{n}}\mathbf{f}({\mathbf{u}}(\tau))d\tau\approx\mathbf{M}^{-1}\Delta\tau_n\mathbf{f}\bigl(\mathcal{F}\bigl(\tau_{n},\tau_{n-1},\mathbf{U}_{n-1}^{(k)},T\bigr)\bigr),
\end{aligned}
\end{equation}
for $n=1,\dots,N.$ In the general case, the system of equation \eqref{kulchytska-ruchka:PPPC_unknownT} is nonlinear and implicit, which requires an additional linearization. 

Building upon the ideas presented in \cite{Kulchytska-Ruchka_2019ag}, which dealt with the time-periodic problem for a known given period $T$, we encorporate an additive splitting of the system matrix in \eqref{kulchytska-ruchka:PPPC_unknownT}. For this let us introduce a modified coarse propagator $\bar{\mathcal{G}},$ which instead of \eqref{kulchytska-ruchka:IVP_window} solves an approximate model with a linearized function $\bar{\mathbf{f}}(\mathbf{u})=\mathbf{A}\mathbf{u}+\mathbf{c}$ on the RHS, i.e.,
\begin{equation}
\label{kulchytska-ruchka:IVP_window_linear}
\begin{aligned}
\mathbf{M}{\mathbf{u}}^{\prime}_n(\tau)&=T\bar{\mathbf{f}}(\mathbf{u}_n(\tau))=T[\mathbf{A}_n\mathbf{u}_n(\tau)+\mathbf{c}_n],\ \ \tau\in(\tau_{n-1},\tau_n]\\
\mathbf{u}_n(\tau_{n-1})&=\mathbf{U}_{n-1}
\end{aligned}
\end{equation} 
with a given {Jacobi-}matrix $\mathbf{A}_n$ and a vector $\mathbf{c}$. Having the linear coarse model we construct a fixed point iteration: for $s=0,1,\dots$ 
\begin{equation}
\label{kulchytska-ruchka:PPPC_unknownT_fxdPt}
\setlength{\arraycolsep}{2.5pt}
\begin{bmatrix}
-\mathbf{I} &    &  & {\bar{\mathcal{G}}_N\left(\cdot,T\right)} & \mathbf{g}^{(k)}_{N} \\
{\bar{\mathcal{G}}_1\left(\cdot,T\right)}  & -\mathbf{I} &  &  & \mathbf{g}^{(k)}_{1}\\
 &  \ddots   &  \ddots & &\vdots \\
&     &  {\bar{\mathcal{G}}_{N-1}\left(\cdot,T\right)} & -\mathbf{I} & \mathbf{g}^{(k)}_{N-1}
\end{bmatrix}
	\begin{bmatrix}
	\mathbf{U}_0^{(k+1,s+1)}\\
	\mathbf{U}_1^{(k+1,s+1)}\\
	\vdots\\
	\mathbf{U}_{N-1}^{(k+1,s+1)}\\
	T^{(k+1)}
	\end{bmatrix}=
	\begin{bmatrix}
		\mathbf{h}_N^{(k+1,s)}\\
		\mathbf{h}_1^{(k+1,s)}\\
		\vdots\\
		\mathbf{h}_{N-1}^{(k+1,s)}
\end{bmatrix}
\end{equation}
where we denote $\mathbf{h}_n^{(k+1,s)}:=\mathbf{b}_n^{(k)}+\bar{\mathcal{G}}\bigl(\tau_n,\tau_{n-1},\mathbf{U}^{(k+1,s)}_{n-1},T\bigr)-{\mathcal{G}}\bigl(\tau_n,\tau_{n-1},\mathbf{U}^{(k+1,s)}_{n-1},T\bigr)$ and {$\bar{\mathcal{G}}_n\bigl(\cdot,T\bigr):=\bar{\mathcal{G}}\bigl(\tau_n,\tau_{n-1},\cdot,T\bigr)$} for $n=1,\dots,N.$ Assuming that $\bar{\mathcal{G}}$ solves \eqref{kulchytska-ruchka:IVP_window_linear} with the implicit Euler method using a single step on $(\tau_{n-1},\tau_n]$ and that all the windows have the same length $\Delta\tau$, we have an explicit representation for the coarse solution 
\begin{equation}
\label{kulchytska-ruchka:Euler_linear}
\left[1/\Delta\tau\cdot\mathbf{M}-T^{(k)}\mathbf{A}\right]\bar{\mathcal{G}}\bigl(\tau_n,\tau_{n-1},\mathbf{U}^{(k+1,s)}_{n-1},T\bigr)=1/\Delta\tau\cdot\mathbf{M}\mathbf{U}^{(k+1,s)}_{n-1}+T^{(k)}\mathbf{c}_n,
\end{equation}
for $n=1,\dots,N.$ Denoting by $\mathbf{C}:=1/\Delta\tau\cdot\mathbf{M}$ and $\mathbf{Q}^{(k)}:=\mathbf{C}-T^{(k)}\mathbf{A}$ 
and plugging this into the system \eqref{kulchytska-ruchka:IVP_window_linear} we obtain
\begin{equation*}
\label{kulchytska-ruchka:PPPC_unknownT_fxdPt_explicit}
\setlength{\arraycolsep}{2.5pt}
\begin{bmatrix}
-\mathbf{Q}^{(k)} &  & & \mathbf{C} & \mathbf{Q}^{(k)}\mathbf{g}^k_{N} \\
\mathbf{C}  & -\mathbf{Q}^{(k)} &  &  & \mathbf{Q}^{(k)}\mathbf{g}^k_{1}\\
 &  \ddots   &  \ddots  & & \vdots\\
 &     &  \mathbf{C} & -\mathbf{Q}^{(k)} & \mathbf{Q}^{(k)}\mathbf{g}^k_{N-1}
\end{bmatrix}
	\begin{bmatrix}
	\mathbf{U}_0^{(k+1,s+1)}\\
	\mathbf{U}_1^{(k+1,s+1)}\\
	\vdots\\
	\mathbf{U}_{N-1}^{(k+1,s+1)}\\
	T^{(k+1)}
	\end{bmatrix}=
	\begin{bmatrix}
		\mathbf{Q}^{(k)}\mathbf{h}_N^{(k+1,s)}-T^{(k)}\mathbf{c}_N\\
		\mathbf{Q}^{(k)}\mathbf{h}_1^{(k+1,s)}-T^{(k)}\mathbf{c}_1\\
		\vdots\\
		\mathbf{Q}^{(k)}\mathbf{h}_{N-1}^{(k+1,s)}-T^{(k)}\mathbf{c}_{N-1}
\end{bmatrix}.
\end{equation*}
\begin{rem}
We note that when the period $T$ is given within the problem setting \eqref{kulchytska-ruchka:TP_pbm}, 
the corresponding block-cyclic matrix (system matrix of \eqref{kulchytska-ruchka:PPPC_unknownT_fxdPt_explicit} without the last column) can be transformed into a block-diagonal using the frequency domain transformation \cite{Biro_2006aa}. A detailed description of the approach as well as a Newton-like linearization of the periodic system within the parallel-in-time setting is presented in \cite{Kulchytska-Ruchka_2019ag}.
\end{rem}

\section{Numerical example}
\label{kulchytska-ruchka:numerics}
\begin{figure}[t]
\begin{minipage}{0.6\textwidth}
  \centering
  \includegraphics[width=0.98\textwidth]{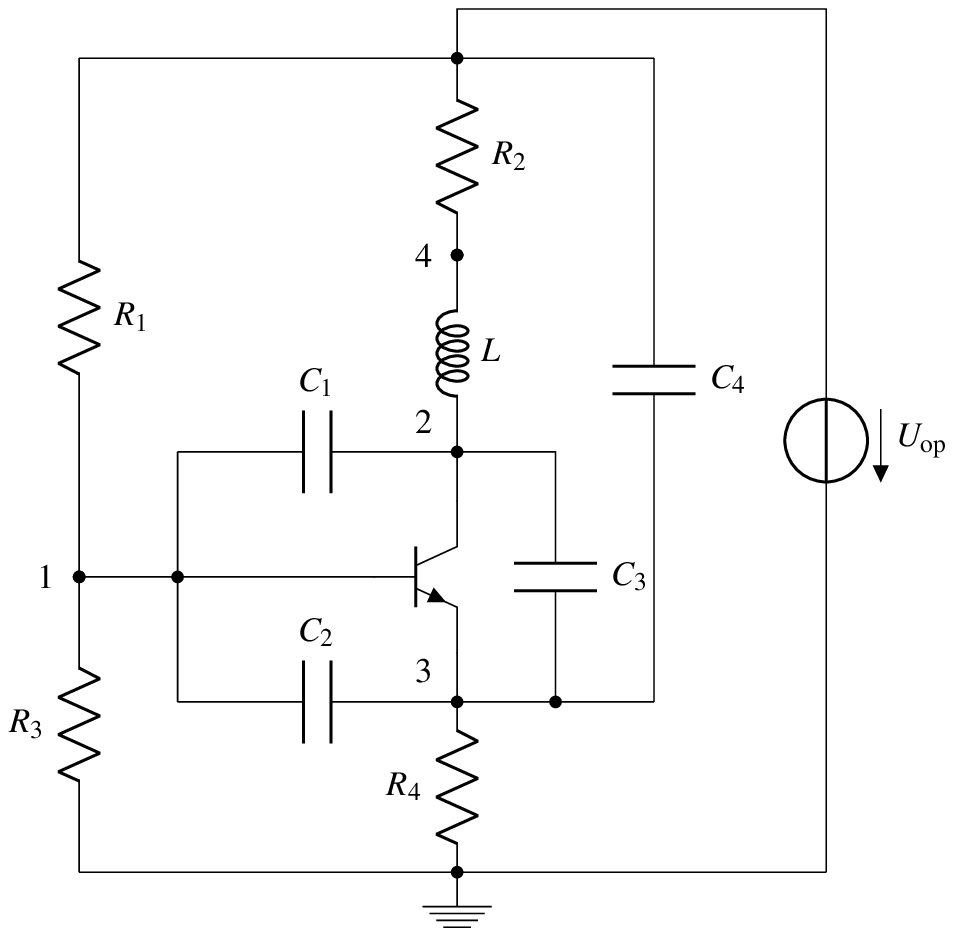}  
\end{minipage}
\begin{minipage}{0.4\textwidth}
\begin{equation*}
		\setlength{\tabcolsep}{12pt}
		\begin{tabular}{l l}
		$C_1=\SI{50}{\pico\farad},$ & $C_2=\SI{1}{\nano\farad},$\\
		$C_3=\SI{50}{\nano\farad},$ & $C_4=\SI{100}{\nano\farad},$\\
		 & \\
		$R_1=\SI{12}{\kilo\ohm},$ & $R_2=\SI{3}{\ohm},$\\
		$R_3=\SI{8.2}{\kilo\ohm},$ & $R_4=\SI{1.5}{\kilo\ohm},$\\
		 & \\
		$L=\SI{10}{\milli\henry},$ & $U_{\mathrm{op}}=\SI{10}{\volt}.$
		\end{tabular}
\end{equation*}
\end{minipage}
\caption{Circuit of the Colpitts oscillator model \cite{Kampowsky_1992aa}.}
\label{kulchytska-ruchka:fig_2} 
\end{figure}

We now consider the Colpitts oscillator model presented in \cite{Kampowsky_1992aa}. It is described by the circuit illustrated in Fig.~\ref{kulchytska-ruchka:fig_2}, which consists of an inductance, a bipolar transistor, as well as of four capacitances and four resistances. The Colpitts oscillator model was exploited in the multi-rate context in \cite{Pulch_2005aa}. 

The mathematical model of the circuit is given by an implicit system of ODEs \cite{Kampowsky_1992aa}, namely, we search for the four node voltages $\mathbf{U}=[U_1,U_2,U_3,U_4]^{\top}$ s.t.
\begin{equation}
\begin{aligned}
&\begin{bmatrix}
1 & 0 & 0 & 0\\
0 & C_1+C_3 & -C_3 & -C_1\\
0 & -C_3 & C_2+C_3+C_4 & -C_2\\
0 & -C_1 & -C_2 & C_1+C_2
\end{bmatrix}
\begin{bmatrix}
\dot{U}_1\\
\dot{U}_2\\
\dot{U}_3\\
\dot{U}_4
\end{bmatrix}\\
&=
\begin{bmatrix}
(U_2-U_1)R_2/L\\ 
(U_{\mathrm{op}}-U_1)/R_2+x_{\mathrm{C}}h(U_4-U_2)-I_{\mathrm{S}}h(U_4-U_3)\\
-U_3/R_4+x_{\mathrm{E}}h(U_4-U_3)-I_{\mathrm{S}}h(U_4-U_2)\\
-U_4/R_3+(U_{\mathrm{op}}-U_4)/R_1-y_{\mathrm{E}}h(U_4-U_3)-y_{\mathrm{C}}h(U_4-U_2)
\end{bmatrix},
\end{aligned}
\end{equation}
with the parameters	
		$y_{\mathrm{E}}=\SI{10}{\micro\ampere},$ $x_{\mathrm{E}}=\SI{1.01}{\milli\ampere},$ 
		{$I_{\mathrm{S}}=\SI{1}{\milli\ampere}$}
		$y_{\mathrm{C}}=\SI{20}{\micro\ampere},$  $x_{\mathrm{C}}=\SI{1.02}{\milli\ampere},$ 
and the nonlinear function $h(x)=\exp(x/U_{\mathrm{T}})-1$, $U_{\mathrm{T}}=\SI{2.585}{\volt},$ coming from the applied transistor model. 
Compared to the model introduced in \cite{Kampowsky_1992aa}, the value of $U_{\mathrm{T}}$ is chosen bigger to {ease the convergence of PP-PC} using the function $h.$ {In practice, one may need appropriate homotopy or damping strategies, see \cite{Deuflhard_2004aa}.}
\begin{figure}[t]
\centering
	\mbox{\includegraphics[width=0.47\linewidth]{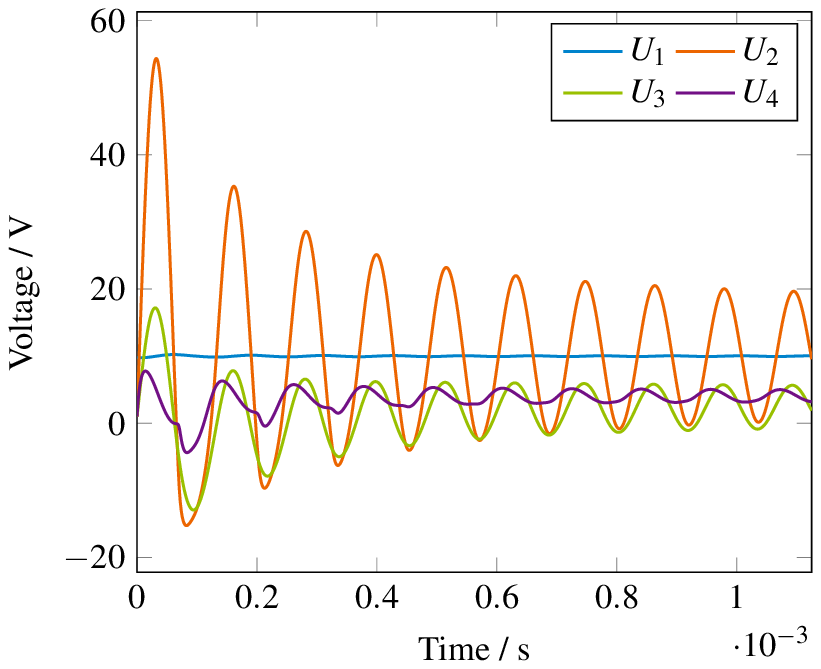}
			\hspace*{0.03\linewidth}
		\includegraphics[width=0.475\linewidth]{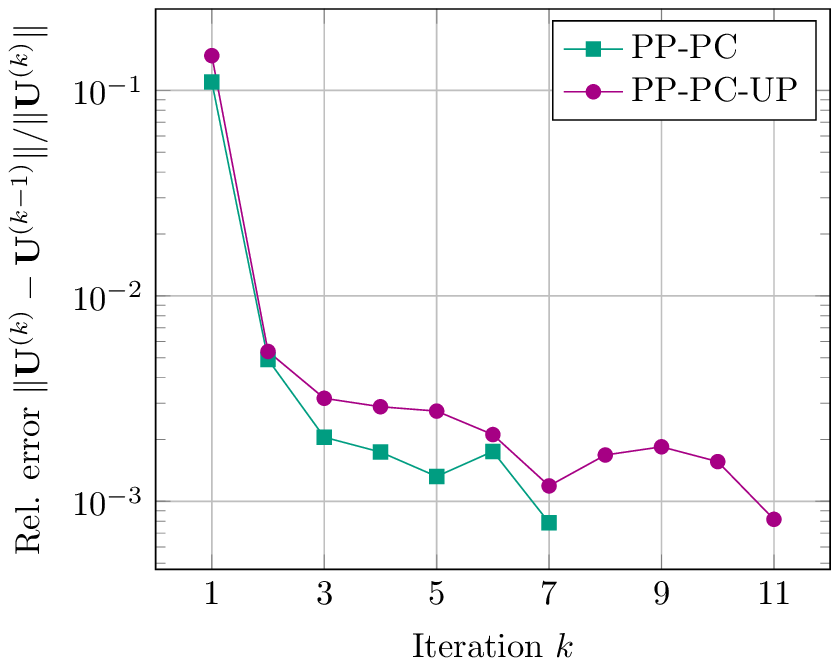}
		}
\caption{Left: Transient behavior of the Colpitts oscillator until the steady state. Right: Convergence of the PP-PC approach with a linearized coarse grid problem for the case when the period $T$ is given \cite{Kulchytska-Ruchka_2019ag} and of PP-PC-UP when $T$ is unknown \eqref{kulchytska-ruchka:PPPC_unknownT_fxdPt_explicit}.}
\label{kulchytska-ruchka:fig_34}   
\end{figure}
The transient behavior of the oscillator on $[0,1.125]$ \SI{}{\milli\second} is shown in Fig.~\ref{kulchytska-ruchka:fig_34} on the left. The time step $\delta T=\SI{0.1125}{\micro\second}$ and the initial value at $t=0$ is $\mathbf{u}_0=[9.75,1,1,1]^{\top}$ are chosen.

To find the periodic steady-state solution and the corresponding period $T$ we apply the iteration \eqref{kulchytska-ruchka:PPPC_unknownT_fxdPt}. Linearization of the nonlinear periodic system on the coarse level is performed using a surrogate linear model, i.e., $\bar{\mathcal{G}}$ solves the problem \eqref{kulchytska-ruchka:IVP_window_linear} with $\bar{\mathbf{f}}(\mathbf{U})=\mathbf{A}\mathbf{U}+\mathbf{c}$ given by
\begin{align}
\mathbf{A}&=
\begin{bmatrix}
-R_2/L & R2/L & & \\
-1/R_2 & -x_{\mathrm{C}}/\bar{U}_{\mathrm{T}} & I_{\mathrm{S}}/\bar{U}_{\mathrm{T}} & (x_{\mathrm{C}}-I_{\mathrm{S}})/\bar{U}_{\mathrm{T}}\\
0 & I_{\mathrm{S}}/\bar{U}_{\mathrm{T}} & -1/R_4-x_{\mathrm{E}}/\bar{U}_{\mathrm{T}} & (x_{\mathrm{E}}-I_{\mathrm{S}})/\bar{U}_{\mathrm{T}}\\
0 & y_{\mathrm{C}}/\bar{U}_{\mathrm{T}} & y_{\mathrm{E}}/\bar{U}_{\mathrm{T}} & -1/R_3-1/R_1-y_{\mathrm{E}}/\bar{U}_{\mathrm{T}}-y_{\mathrm{C}}/\bar{U}_{\mathrm{T}}
\end{bmatrix},\\
\mathbf{c}&=\left[0,U_{\mathrm{op}}/R_2,0,U_{\mathrm{op}}/R_1\right]^{\top},
\end{align}
with $\bar{U}_{\mathrm{T}}=\SI{0.2585}{\volt}.$ The unit interval $[0,1]$ is split into $N=10$ windows. The coarse time step $\Delta\tau=0.1$ and the fine step $\delta\tau=10^{-4}$ were chosen within the time integration. The calculated period with the fixed point iteration \eqref{kulchytska-ruchka:PPPC_unknownT_fxdPt_explicit} is $T=\SI{0.1125}{\milli\second}.$ The right-hand part of the Fig.~\ref{kulchytska-ruchka:fig_34} shows convergence of the PP-PC iteration with the linearization from \cite{Kulchytska-Ruchka_2019ag} for a given period $T$ as well as 
for an unknown period {(PP-PC-UP)}. 
Both results are obtained up to the relative tolerance of $10^{-3}.$ One can see that in case when $T$ is known the method required less iterations, as {one would} expect. {When comparing the computational cost of the computations in terms of the number of linear systems solves, PP-PC and PP-PC-UP delivered the periodic solution effectively $4$ and $3$ times faster than the sequential time stepping, respectively. }

\section{Conclusions}
\label{kulchytska-ruchka:conclusions}
An iterative parallel-in-time method for solving time-periodic 
problems where the period is not initially given 
is proposed in this paper. It complements   
the system of equations originating from the Parareal-like algorithm 
with an additional variable $T$ and gives an underdetermined system of nonlinear equations. A linearization using the fixed point iteration is applied on the coarse grid. The algorithm is verified through its application to a Colpitts oscillator model.

\begin{acknowledgement}
The authors thank Roland Pulch from Universit\"at Greifswald for his assistance with implementation and for the fruitful discussions on the Colpitts oscillator model.

This research was supported by the Excellence Initiative of the German 
Federal and State Governments and the Graduate School of Computational 
Engineering at Technische Universit\"at Darmstadt, as well as by DFG 
grant SCHO1562/1-2 and BMBF grant 05M2018RDA (PASIROM). 
\end{acknowledgement}

\bibliographystyle{plain}



\end{document}